\documentclass{proc-l}

\usepackage{amssymb}
\usepackage{amsmath}
\usepackage{amsthm}
\usepackage{amsbsy}
\usepackage[dvips]{graphics}
\usepackage{amscd}
\usepackage{bm}

\theoremstyle{plain}
\newtheorem{theorem}{Theorem}[section]

\newtheorem{lemma}[theorem]{Lemma}
\newtheorem{corollary}[theorem]{Corollary}
\newtheorem{proposition}[theorem]{Proposition}

\theoremstyle{definition}
\newtheorem{definition}[theorem]{Definition}
\newtheorem{example}[theorem]{Example}

\theoremstyle{remark}
\newtheorem{remark}[theorem]{Remark}

\newcommand{\del}{\partial}

\newcommand{\F}{\mathfrak{F}}
\newcommand{\Z}{{\mathbb{Z}}}
\newcommand{\C}{{\mathbb{C}}}
\newcommand{\R}{{\mathbb{R}}}
\renewcommand{\P}{{\mathcal{P}}}
\renewcommand{\H}{{\mathbb{H}}}

\newcommand{\Zz}{\mathbb{Z}/2\mathbb{Z}}

\renewcommand{\i}{\bm{i}}
\renewcommand{\j}{\bm{j}}
\renewcommand{\k}{\bm{k}}

\newcommand{\abs}[1]{\lvert#1\rvert}

\begin{document} 

\title{Contact structures on elliptic $3$-manifolds}

\begin{abstract}
We show that an oriented elliptic $3$-manifold admits a universally
tight positive contact structure if and only if the corresponding
group of deck transformations on $S^3$ (after possibly conjugating by
an isometry) preserves the standard contact structure.

We also relate universally tight contact structures on $3$-manifolds
covered by $S^3$ to the isomorphism $SO(4)=(SU(2)\times SU(2))/{\pm
1}$.

The main tool used is equivariant framings of $3$-manifolds.

\end{abstract}


\date{\today}
\author{Siddhartha Gadgil}

\address{	Department of mathematics\\
		SUNY at Stony Brook\\
		Stony Brook, NY 11794}
\email{gadgil@math.sunysb.edu}

\maketitle 

A contact structure $\xi$ on a $3$-dimensional manifold $M$ is a
smooth, totally non-integrable tangent plane field, i.e., a tangent
plane field $\xi$ locally of the form $\xi=ker(\alpha)$ for a $1$-form
$\alpha$ such that $\alpha\wedge d\alpha$ is everywhere
non-degenerate. We shall assume that $M$ is oriented. We say $\xi$ is
positive if the orientation on $M$ agrees with that induced by the
volume form $\alpha\wedge d\alpha$. Observe that the orientation of
$\alpha\wedge d\alpha$ does not depend on the sign of $\alpha$, and is
thus determined by $\xi$ (even though $\xi=ker(\alpha)$ only locally).

A central role in understanding $3$-dimensional manifolds has been
played by co-dimension one structures -- surfaces, foliations and
laminations -- in these manifolds. Without additional conditions such
structures always exist, and are of not much consequence. However, the
presence of \emph{essential} co-dimension one structures --
\emph{incompressible} surfaces, \emph{taut} foliations and
\emph{essential} laminations, leads to deep topological consequences.

By the work of Eliashberg~\cite{El}\cite{El2}, there is a similar
dichotomy among contact structures between  \emph{tight}
contact structures and  \emph{overtwisted} contact structures
(we recall the definitions in the next section). Further, there are deep
connections between taut foliations and tight contact structures by
the work of Eliashberg and Thurston~\cite{ET}, and more recently of
Honda, Kazez and Matic~\cite{HKM}\cite{HKM2}.

There is however one significant difference between contact structures
and the other co-dimension one structures -- while one of the most basic
consequences of the existence of other essential co-dimension one
structures in $M$ is that the universal cover of $M$ is $\R^3$ (this
was in fact used to demonstrate the utility of essential laminations
when they were introduced), one of the most basic examples of a tight
contact structure is the \emph{standard} contact structure on $S^3$
(which we recall in the next section). Several quotients of $S^3$ also
admit tight contact structures.

Thus, tight contact structures clearly reveal a different aspect of
$3$-manifolds, at least in some cases. Our main goal here is to relate
the existence of contact structures on elliptic manifolds (i.e.,
quotients of $S^3$ by a group of isometries) with tight universal
covers to the  isomorphism $SO(4)=(SU(2)\times SU(2))/{\pm
1}$, and more generally to spherical structures.

Our first main result is the following.

\begin{theorem}\label{T:main}
Suppose $M=S^3/G$ where $G$ is a group of isometries of $S^3$. Then
the oriented manifold $M$ has a positive contact structure with tight
universal cover if and only if $G$ (after possibly conjugating by an
isometry) leaves invariant the standard contact structure on $S^3$.
\end{theorem}

It is easy to deduce, using the isomorphism, which finite groups
$G\subset SO(4)$ preserve the standard contact structure on $S^3$. The
main content of this paper is to show that the other elliptic
manifolds do not admit a positive contact structure with tight
universal cover. The main tool we use is \emph{equivariant framings},
introduced in~\cite{Ga}, and the methods of this paper are essentially
a straightforward extension of the ones in that paper. An alternative
approach to the results we obtain is to use Gompf's invariants for
tangent plane fields (see~\cite{Go}). We refer to~\cite{Ga} for the
relation between these and equivariant framings.

The same methods yield the following stronger result.

\begin{theorem}\label{T:main2}
Suppose $M$ is a $3$-manifold with universal cover $S^3$ and suppose
$\pi_1(M)\cong G$ for some $G\subset SO(4)$ that acts freely on
$S^3\subset \R^4$, $G$ not cyclic. Then if there are positive contact
structures with tight universal cover for both orientations of $M$,
then there are positive contact structures with tight universal cover
for both orientations of $S^3/G$. Moreover $G$ preserves the standard
contact structure on $S^3$ (after possibly conjugating by an
isometry).
\end{theorem}

Thus the restrictions to existence of positive contact structures with
tight universal cover on elliptic manifolds are essentially homotopy
theoretic. 

\begin{remark}
An elliptic $3$-manifold $N=S^3/G$ with $G$ not cyclic is determined
up to (not necessarily orientation preserving) isometry by its
fundamental group. Thus, the above result says that if $M$ admits
contact structures for both orientations, so does the unique elliptic
$N$ with the same fundamental group.
\end{remark}

From another point of view, we conclude that if $M=S^3/G$ is an
elliptic manifold, with $G\subset SO(4)=(SU(2)\times SU(2))/{\pm 1}$
not cyclic, then the image of $G$ under projection onto the two
factors is a \emph{topological} property of $M$ determined by contact
structures on $M$.

There is an even closer connection between the  isomorphism
and \emph{co-orientable} positive contact structures with tight
universal cover whose Euler class is trivial on the quotient manifold
as we see in the following theorem.

\begin{theorem}\label{T:left}
Suppose $M=S^3/G$ is an elliptic manifold and $G$ is not a cyclic
group. Then $M$ admits a co-orientable positive contact structure
with tight universal cover whose Euler class is trivial if and only if
$G\subset SU(2)\times 1$.
\end{theorem}

\section{Preliminaries}

We recall below basic definitions and results in contact geometry that
we need. For details and motivations we refer to~\cite{El}. Henceforth
let $M$ denote a closed oriented $3$-manifold.

\begin{definition}
A contact structure $\xi$ on $M$ is a smooth tangent plane field that
is locally of the form $\xi=ker(\alpha)$ for a $1$-form $\alpha$ such
that $\alpha\wedge d\alpha$ is everywhere non-degenerate. We say that
$\xi$ is positive if the orientation on $M$ agrees with that induced
by $\alpha\wedge d\alpha$.
\end{definition}

An important example is the \emph{standard contact structure} on
$S^3$.

\begin{example}
Consider the tangent plane field $\xi$ on $S^3\subset{\C^2}$ that is
perpendicular to the vector field $V:S^3\to TS^3$, $V:(z_1,z_2)\mapsto
(\i z_1, \i z_2)$. This is a positive contact structure called the
\emph{standard} contact structure on $S^3$.
\end{example}

\begin{definition}\label{D:overtwist}
A contact structure $\xi$ on $M$ is said to be \emph{overtwisted} if
there is an embedded disc $D\subset M$ so that $TD|_{\del
D}=\xi|_{\del D}$. A contact structure that is not overtwisted is said
to be \emph{tight}.
\end{definition}

\begin{definition}
A contact structure $\xi$ on $M$ is said to be \emph{universally
tight} if the pullback of $\xi$ to the universal cover of $M$ is
tight.
\end{definition}

Universally tight contact structures are tight as any cover of an
overtwisted contact structure is overtwisted. This follows as the disc
$D$ in Definition~\ref{D:overtwist} lifts to any cover.

Eliashberg~\cite{El2} has shown that there is a unique overtwisted
contact structure representing each homotopy class of tangent plane
fields on a manifold $M$. This is far from true in the case of tight
contact structures, and their existence is still mysterious.

The standard contact structure on $S^3$ is tight. This is essentially
the only tight contact structure on $S^3$ by the following result of
Eliashberg~\cite{El2}.

\begin{theorem}[Eliashberg]
Any positive tight contact structure on $S^3$ is isotopic to the
standard one.
\end{theorem}

\section{Elliptic $3$-manifolds}

An elegant classification of elliptic $3$-manifolds is obtained by
Hopf~\cite{Ho} (see Scott~\cite{Sc} for a very readable account) using
the isomorphisms $SO(4)=(SU(2)\times SU(2))/{\pm 1}$ and
$SO(3)=SU(2)/{\pm 1}$. We outline this in this section. This has a
transparent connection with contact geometry which we shall exploit to
construct contact structures. For proofs we refer to~\cite{Sc}.

\subsection{The exceptional isomorphisms}

Consider the quaternions
$\H=\{\omega_1+\j\omega_2:\omega_1,\omega_2\in\C\}$. The group
$SU(2)=S^3$ can be identified with the set of unit quaternions. This
acts on the quaternions isometrically by left multiplication and by
right multiplication, giving a surjective map $SU(2)\times SU(2)\to
SO(4)$ with kernel $\{1,-1\}$. This gives the  isomorphism
$\phi:(SU(2)\times SU(2))/{\pm 1}\overset{\cong}\to SO(4)$. The
isomorphism $SO(3)=SU(2)/{\pm 1}$ is obtained by considering the
action of $S^3$ on itself by conjugation.

Notice that $\H$ has a complex structure induced by right
multiplication by $\C\subset\H$. The image of $S^3\times S^1$ gives
complex linear maps as these commute with every element of
$1\times\C$.

\subsection{The classification}

An elliptic manifold is a quotient of $S^3$ by a finite subgroup
$G\subset SO(4)$ that acts without fixed points on the unit
sphere. The above isomorphisms can be used together with the
classification of finite subgroups of $SO(3)$ to classify such groups.

Recall that finite subgroups of $SO(3)$ are the cyclic groups,
dihedral groups and the groups of symmetries of the tetrahedron,
octahedron and icosahedron. The inverses images of these groups under
the covering $SU(2)\to SO(3)$ give cyclic groups, quaternionic groups
$Q_{4n}$, the binary tetrahedral group $T$, the binary octahedral
group $O$ and the binary icosahedral group $I$ respectively.

Using the fact that $G$ acts freely, one can deduce the following
proposition (this is Theorem 4.10 of~\cite{Sc}) which is important for
our purposes.

\begin{proposition}
$G$ is conjugate to the image in $SO(4)$ of a subgroup $\tilde G$ of
$S^1\times S^3$ or of a subgroup of $S^3\times S^1$. 
\end{proposition}

Suppose now that $G$ is a subgroup of $\phi(S^1\times S^3)$ that acts
freely on $S^3$. We have the following classification.

\begin{theorem}[Hopf]\label{T:class}
Suppose $G\subset \phi(S^1\times S^3)$ acts freely on $S^3$. Then one
of the following holds.
\begin{enumerate}

\item\label{Cy} $G$ is cyclic

\item\label{Qt} $G$ is the product of a quaternionic group $Q_{4n}$ in
$\phi(1\times S^3)$ with a cyclic group of relatively prime order in
$\phi(S^1\times 1)$.

\item\label{T} $G$ is the product of $T\subset \phi(1\times S^3)$ with
a cyclic group of relatively prime order in $\phi(S^1\times 1)$

\item\label{O} $G$ is the product of $O\subset \phi(1\times S^3)$ with
a cyclic group of relatively prime order in $\phi(S^1\times 1)$

\item\label{I} $G$ is the product of $I\subset \phi(1\times S^3)$ with
a cyclic group of relatively prime order in $\phi(S^1\times 1)$

\item\label{Q2} $G$ is the quotient under $\phi$ of an index $2$
diagonal subgroup (for the definition of diagonal subgroups see Scott,
page 453) of $C_{2m}\times Q_{4n}$, where $m$ is odd and $n$ and $m$
are relatively prime.

\item\label{T2} $G$ is an index $3$ diagonal subgroup of $C_{3m}\times
T$, where $m$ is odd.

\end{enumerate}

\end{theorem}

\subsection{Reversing orientations}

Suppose $M$ is an oriented elliptic manifold of the form $S^3/G$. Then
it follows using, for instance, the orientation reversing
anti-isomorphism
	$$\psi:\omega_1+\j\omega_2\mapsto\omega_1+\j\bar{\omega_2}$$
that the oriented manifold $-M$ ($M$ with the opposite orientation) is
obtained by switching the left and right components of $G\subset
S^3\times S^3$

\section{Quotient contact structures}

We have seen that the standard contact structure $\xi$ on $S^3\subset
\C^2$ is characterised by being perpendicular to the vector field
$V:(z_1,z_2)\mapsto (\i z_1,\i z_2)$. We can immediately deduce that
several elliptic manifolds have quotient contact structures.

\begin{lemma}
Suppose $M=S^3/G$ with $G\subset \phi(S^3\times S^1)$, then $M$ has a
co-orientable quotient contact structure induced by $\xi$.
\end{lemma}
\begin{proof}
Any $g\in S^3\times S^1$ acts by complex linear maps and hence
preserves $V$.
\end{proof}

The following is immediate.
\begin{theorem}
Any elliptic manifold $M$ admits a universally tight contact structure
for at least one orientation of $M$.
\end{theorem}
\begin{proof}
Let $M=S^3/G$. If $G\subset \phi(S^3\times S^1)$ then the above lemma
shows that $M$ has such a contact structure. Otherwise $G\subset
\phi(S^1\times S^3)$ and hence the manifold $-M$ obtained by reversing
the orientation on $M$ is a quotient of $S^3$ by a subgroup of
$S^3\times S^1$ and hence has a universally tight contact structure.
\end{proof}

Allowing deck transformations that reverse the co-orientation, we can
construct more quotient contact structures.

\begin{lemma}
Suppose $M=S^3/G$ with $G\subset \phi(S^3\times (S^1\cup \j S^1))$,
then $M$ has a quotient contact structure induced by $\xi$.
\end{lemma}
\begin{proof}
Any $g\in S^3\times (S^1\cup \j S^1)$ acts by complex linear or
anti-linear maps and hence either preserves $V$ or takes $V$ to
$-V$. In either case $\xi$ is preserved.
\end{proof}

\begin{theorem}
The manifolds in cases \ref{Cy}, \ref{Qt} and \ref{Q2} of
Theorem~\ref{T:class} admit universally tight contact structures in
each orientation.
\end{theorem}
\begin{proof}
Let $M$ be such a manifold. In case~\ref{Cy}, after conjugation
$\pi_1(M)\subset \phi(S^1\times S^1$). 

Next, let $\zeta_n\in\C$ denotes a primitive $n$th root of unity. In
case~\ref{Qt}, $\pi_1(M)$ is the subgroup of $SO(4)$ that is the image
of $(\j,1)\in S^3\times S^3$ and $(\zeta_{2n}, 1)\in S^3\times S^3$, or
the product of this with the subgroup generated by $(1,\zeta_m)$ for
some $m$. In case~\ref{Q2}, $\pi_1(M)$ is the subgroup of $SO(4)$ that
is the image of $(\j,\zeta_{2k})\in S^3\times S^3$ and $(\zeta_{2n}, 1)\in
S^3\times S^3$, or the product of this with the subgroup generated by
$(1,\zeta_m)$ for some $m$.

Thus for one orientation $\pi_1(M)\subset \phi(S^3\times S^1)$ and for the
other $\pi_1(M)\subset \phi(S^1\times (S^1\cup \j S^1))$. In either case we
get a quotient contact structure.
\end{proof}

The main content of this paper is that there is no universally tight
contact structure for the remaining elliptic manifolds.

\section{Equivariant framings}

The main tool we use is the \emph{equivariant framing} of a
$3$-manifold of the form $S^3/G$ introduced in~\cite{Ga}. We outline
below the relevant concepts and results.

We use the fact that the tangent bundle of an oriented $3$-manifold is
trivialisable. Let $M=S^3/G$, where $G$ is a finite group acting
without fixed points on $S^3$ (not necessarily by isometries).

We define an invariant $\F(M)$ of $M$ with a given orientation, which
we call the \emph{equivariant framing} of $M$. Recall that the
homotopy classes of trivialisations of the tangent bundle of $S^3$ are
a torseur of $\Z$ (i.e., a set on which $\Z$ acts freely and
transitively), which can moreover be canonically identified with $\Z$
by using the Lie group structure of $S^3$ as the unit quaternions and
identifying a left-invariant framing with $0\in\Z$.

Now, find a trivialisation $\tau$ of $TM$ and pull it back to one of
$TS^3$. Under the above identification, this gives an element
$\F(M,\tau)\in\Z$. This depends on $\tau$, but we can see that its
reduction modulo $\abs{G}$, when $H_1(M,\Zz)=0$ (in particular when
$\abs{G}$ is odd), and modulo $\abs{G}/2$ otherwise, is well-defined
by the following straightforward proposition.

\begin{proposition}[see~\cite{Ga}] 
Suppose $\tau_i,i=1,2$ are trivialisations of $TM$ and $\pi^*(\tau_i)$
are their pullbacks under the covering map $\pi:S^3\to M$.  Then
$\pi^*(\tau_1)-\pi^*(\tau_2)$ is divisible by $\abs{G}$ when
$H_1(M,\Zz)=0$ and by $\abs{G}/2$ when $H_1(M,\Zz)\neq 0$.
\end{proposition}
\qed

Using this, we define the invariant $\F(M)$.

\begin{definition} Let $M=S^3/G$, where $G$ is a finite
group acting without fixed points on $S^3$. The equivariant framing
$\F(M)\in \Z/\left<G\right>\Z$, where $\left<G\right>=\abs{G}$ when
$H_1(M,\Zz)=0$ and $\left<G\right>=\abs{G}/2$ when $H_1(M,\Zz)\neq 0$, is
the equivalence class of the trivialisation of $TS^3$ obtained by
pulling back a trivialisation of $TM$.
\end{definition}

The above definition does not depend on the identification of $S^3$
with the universal cover of $M$, since two such identifications differ
by an orientation preserving self-homeomorphism of $S^3$, which is
isotopic to the identity.

We shall need the following results regarding the framing invariant.

\begin{proposition}
Suppose $M=S^3/G$. If $G\subset \phi(S^3\times 1)$ (respectively $G\subset
\phi(1\times S^3)$) then $\F(M)=0$ (respectively $\F(M)=1$)
\end{proposition}
\begin{proof}
If $G\subset \phi(S^3\times 1)$, as the left-invariant trivialisation is
preserved by the action of $G$, it gives a framing $\tau$ on the
quotient $M$. The pullback of $\tau$ is then the left-invariant
framing, i.e., the $0$ element.

On the other hand, the right-invariant trivialisation is preserved by
the action of a group $G\subset \phi(1\times S^3)$, and hence gives a
framing $\tau'$ on $M=S^3/G$ that pulls back to the right-invariant
trivialisation on $M$. 

Thus, $\F(M)$ is determined by the map $\psi:S^3\to SO(3)$ that at a
point $g\in S^3$ is the matrix of transition between the ordered basis
$(\hat u,\hat v,\hat w)=(g\i,g\j,g\k)$ of $T_gS^3$ corresponding to
the left-invariant trivialisation, and the ordered basis $(\hat
u',\hat v',\hat w')=(\i g,\j g,\k g)$ corresponding to the
right-invariant trivialisation. More precisely, $F(M)$ is the degree
of the lift of this map to the \emph{comparison map}
$\tilde\psi:S^3\to S^3$.

Observe that $(\hat u',\hat v',\hat w')=g^{-1}(\hat u,\hat v,\hat
w)g$, and hence the map $\phi:S^3\to SO(3)$ at the point $g$ is the
action by conjugation of $g$ on the Lie algebra of $S^3$. Thus, $\psi$
is the standard $2$-fold covering map $\phi:S^3\to SO(3)$ given by the
adjoint action. It follows that $\tilde\psi$ is the identity map and
hence has degree one.
\end{proof}

\begin{corollary}\label{T:lftrt}
If $M=S^3/G$ with $G\subset \phi(S^3\times 1)$ and $\abs{G}>2$ then no
trivialisation of $-M$ pulls back to one isotopic to the left
invariant trivialisation of $S^3$.
\end{corollary}

\section{Non-existence of contact structures}

Suppose now that $M=S^3/G$ and $\xi$ is a positive contact structure
on $M$. We shall associate framings of $M$ to certain contact
structures. Note that the following proposition does not require $G$
to act by isometries.

\begin{proposition}\label{T:ctcfr}
Let $M=S^3/G$ be a manifold with a positive, co-orientable contact
structure $\xi$ with trivial Euler class. Then there is a framing
canonically associated to $\xi$. Further, the pullback of this framing
to any manifold that covers $M$ is the framing induced by the pullback
of the contact structure.
\end{proposition}
\begin{proof} 
Choose and fix a co-orientation for $\xi$. This induces an orientation
on the plane-bundle given by the contact structure, which we identify
with $\xi$.

As the Euler class of $\xi$ is trivial, there is a trivialisation of
$\xi$ as a vector bundle. Further, two trivialisations differ by a map
$M\to S^1$. As $H^1(M)=0$ (as $H_1(M)=G/{[G,G]}$ is finite), any such
map is homotopic to a constant map.

Thus, there is a trivialisation $(X_1,X_2)$ of $\xi$, canonical up to
homotopy. This, together with a vector $X_3$ normal to $\xi$ that is
consistent with the co-orientation, gives a framing
$(X_1,X_2,X_3)$ of $TM$.

The homotopy class of this trivialisation does not depend on the
choice of co-orientation since $(X_1,-X_2,-X_3)$ gives a
trivialisation corresponding to the opposite co-orientation, and this
is clearly homotopic to the trivialisation $(X_1,X_2,X_3)$.

As the trivialisation of $\xi$ pulls back to give a trivialisation of
the pullback to any cover of $\xi$, and the pullback of $\xi$ is
co-orientable and has trivial Euler class, the second claim follows.
\end{proof}

We can now prove our first non-existence result. Let $\P=S^3/I$ be the
Poincare homology sphere with a fixed orientation, and let $-\P$
denote the same manifold with the opposite orientation. The Poincar\'e
homology sphere is the quotient of $S^3$ by a group acting by left
multiplication, and $-\P$ is the quotient of an action by right
multiplication. We can now prove the following theorem.

\begin{theorem}[Gompf, see also \cite{Ga}] 
The manifold $-\P$ does not have a universally tight positive contact
structure.
\end{theorem} 
\begin{proof}
Let $\xi$ be a positive contact structure on $-\P$. As $-\P$ is a
homology sphere, it follows that the first Stiefel-Whitney class of
$\xi$ vanishes, and hence $\xi$ is co-orientable. Thus its Euler class
is a well-defined element of $H^2(M)=0$ and hence must vanish. Thus
the hypotheses of Proposition~\ref{T:ctcfr} are satisfied.

As $-\P$ is the quotient of $S^3$ by a group acting by right
multiplication, it follows that any framing on $-\P$ pulls back to one
homotopic to a framing invariant under right Lie multiplication, or
one differing from this by $\left\abs{\pi_1(\P)\right}$ units (as
$H_1(\P,\Z_2)=0$). However, if $-\P$ had a universally tight positive
contact structure $\xi$, then the pullback of $\xi$ to $S^3$
is isotopic to the standard contact structure . Hence the framing
associated to $\xi$ pulls back to give the framing associated to left
Lie multiplication. This contradicts Corollary~\ref{T:lftrt}
\end{proof}

\begin{remark} 
Etnyre and Honda~\cite{EH} have shown that $-\P$ does not have a
positive tight contact structure (and in particular does not have a
positive universally tight contact structure).
\end{remark}

\begin{remark} 
Gompf has proved the above results using methods that can be seen
to be essentially the same as those of this paper. For details of the
relation to his methods, see~\cite{Ga}.
\end{remark}

We next consider the manifold $M=S^3/G$, where $G\subset \phi(S^1\times
S^3)$ is as in case~\ref{T}, \ref{O} or~\ref{T2} of
Theorem~\ref{T:class}.

\begin{theorem}\label{T:noctc}
Suppose $M=S^3/G$, where $G\subset \phi(S^1\times S^3)$ is as in
case~\ref{T}, \ref{O}, \ref{I} or~\ref{T2} of
theorem~\ref{T:class}. Then $M$ does not admit a universally tight
positive contact structure.
\end{theorem}
\begin{proof}
The result has been proved in the case~\ref{I}. Suppose next that we
are in one of the cases~\ref{T} or~\ref{T2}.

Let $\xi$ be a positive contact structure on $M$. In cases~\ref{T}
and~\ref{T2}, $\pi_1(M)$ is the product of a group with the presentation
$$\langle x,y,z; x^2=(xy)^2=y^2, zxz^{-1}=y, zyz^{-1}=xy,
z^{3^k}=1\rangle$$ 
with a cyclic group of order $n$ for some $n$ that
is odd and not divisible by $3$. By abelianising, we see that
$H_1(M)=\Z/{3^kn\Z}$, which has odd order. Hence $H^1(M,\Zz)=0$. Hence
the first Stiefel-Whitney class of $\xi$ vanishes and $\xi$ is
co-orientable. Thus its Euler class is a well-defined element of
$H^2(M)$.

Let $M'$ be a cover of $M$ corresponding to the subgroup generated by
an element $g\in G$ of order $4$ (such an element exists in these
cases). Then the pullback map $H^2(M)\to H^2(M')=\Z/{4\Z}$ vanishes as
$H^2(M)=H_1(M)$ has odd order. Thus the hypotheses of
Proposition~\ref{T:ctcfr} are satisfied by the pullback $\xi'$ of
$\xi$ to $M'$.

On the other hand, $M'$ is the quotient of $S^3$ by a subgroup acting
by right multiplication, and hence an argument similar to the
icosahedral case gives a contradiction.

Finally a manifold $M$ corresponding to case~\ref{O} has a cover $N$
with fundamental group $T\subset \phi(S^1\times S^3)$. Hence $M$ does
not admit a positive universally tight contact structure as this would
pull back to give a universally tight positive contact structure on
$N$, contradicting the above.

\end{proof}

We are now in a position to complete the proof of our main results.

\begin{proof}[Proof of theorems~\ref{T:main}, \ref{T:main2} and~\ref{T:left}]

We first prove Theorem~\ref{T:main}. Suppose $M=S^3/G$ with $G\subset
SO(4)$ as in the hypothesis. For one of the orientations of $M$
manifold, we have seen that $\pi_1(M)$ fixes the standard contact
structure on $S^3$.

For the other orientation, in the cases \ref{Cy}, \ref{Qt} and
\ref{Q2} of Theorem~\ref{T:class} we have seen that $\pi_1(M)$ fixes
the standard contact structure. On the other hand,
Theorem~\ref{T:noctc} shows that in cases~\ref{T}, \ref{O}, \ref{I}
and~\ref{T2}, the quotient manifold with one of its orientations does
not admit a positive universally tight contact structure. This
exhausts all the cases, proving the result.

We next prove Theorem~\ref{T:left}. Suppose $M=S^3/G$, $G\subset
SU(2)\times 1$. Then $G$ fixes the standard contact structure on $S^3$
and hence induces a universally tight positive contact structure $\xi$
on $M$. Further, the section $g\mapsto g\i$ of the standard contact
structure is fixed by $G$ and hence descends to a section of
$\xi$. Thus $\xi$ has trivial Euler class.

Conversely, suppose $M=S^3/G$ and $G\not\subset SU(2)\times 1$
is not cyclic, and $M$ admits a universally tight positive contact
structure. By considering cases, it follows that some cover of $M$ is
of the form $S^3/H$, $H\subset \phi(1\times SU(2))$, and $\abs{H}>2$. By
the hypothesis that the Euler class of $\xi$ is trivial and
Proposition~\ref{T:ctcfr}, the pullback of $\xi$ to $S^3/H$ has a
framing associated to it. As before, this lifts to a framing isotopic
to the left-invariant framing of $S^3$ implying $\F(S^3/H)=0$, which
gives a contradiction as $H\subset \phi(1\times SU(2))$ and hence
$\F(S^3/H)=1$.

\begin{remark}
In the case of cyclic groups $G$, the above proof still shows that
under the usual hypothesis $\F(M)=0$. For prime order cyclic groups,
by the results of~\cite{Ga} we still can conclude that $G\subset
\phi(S^3\times 1)$. Tight contact structures on lens spaces (i.e., $G$
cyclic) have been completely classified by Giroux~\cite{Gi} and
Honda~\cite{Hon}.
\end{remark}

To prove Theorem~\ref{T:main2} we shall need the following proposition.

\begin{proposition}
Suppose an oriented manifold $M$ of the form $M=S^3/\Gamma$, with
$\Gamma$ a finite group acting freely (but not necessarily
isometrically) on $S^3$, has a framing that lifts to one isotopic to
the left-invariant framing on $S^3$. Then $-M$ has a framing that
lifts to one isotopic to the right-invariant framing on $S^3$.
\end{proposition}
\begin{proof}

Let $(X_1,X_2,X_3)$ be the framing on $M$ that lifts to one isotopic
to the left-invariant framing on $S^3$. Then $(-X_1,-X_2,-X_3)$ is a
framing on $-M$ that lifts to a framing isotopic to a left-invariant
framing on $S^3$ corresponding to the opposite orientation. Here the
universal covering of $-M$ can be naturally identified with $S^3$ with
the opposite orientation.

Composing the universal covering map of $-M$ with an orientation
reversing homeomorphism $\psi:S^3\to S^3$, we get a covering map from
$S^3$ with its usual orientation to $-M$ that respects
orientations. The pullback $\F$ of the framing on $-M$ is isotopic to
the pullback under $\psi$ of a left-invariant framing $(\hat u,\hat
v,\hat w)$ on $S^3$ corresponding to the opposite orientation. As any
two left-invariant framings corresponding to an orientation are
equivalent, we can take $(\hat u,\hat v,\hat w)_g=g(-\i,-\j,-\k)$.

We use the orientation reversing map $\psi:q\to \bar q$, where $q\to
\bar q$ maps $\i\to-\i$, $\j\to-\j$, $\k\to-\k$ and $1\to 1$ and is
linear (over $\R)$ on $\H$. This is an involution and is  anti-linear
over $\H$, i.e., $\overline{pq}=\bar{q}\bar{p}$.

Observe that $(\hat u,\hat v,\hat
w)_g=g(\bar{\i},\bar{\j},\bar{\k})$. Thus, as $\psi:q\to \bar q$ is an
anti-linear involution, the pullback of $(\hat u,\hat v,\hat
w)_g=g(\bar{\i},\bar{\j},\bar{\k})$ under $\psi$ at $\psi^{-1}(g)=\bar
g$ is $(\i \bar g, \j\bar g,\k\bar g)=(\i,\j,\k)\bar g$. This is
precisely the right-invariant framing at this point corresponding to
the usual orientation.

\end{proof}

We can now prove Theorem~\ref{T:main2}. Suppose $M=S^3/\Gamma$, with
$\Gamma$ a finite group acting freely (but not necessarily
isometrically) on $S^3$. By hypothesis, $\Gamma\cong G$ as groups for
some subgroup $G\subset SO(4)$, $G$ not cyclic, that acts freely on
$S^3$.

Suppose $M$ admits a universally tight positive contact structure for
each of its orientations. We shall show that $G$ preserves the
standard contact structure on $S^3$, and hence in particular $S^3/G$
has a universally tight positive contact for each orientation.  By the
above, this is equivalent to showing that $G$ is not in one of the
cases~\ref{T}, \ref{O}, \ref{I} and~\ref{T2}.

We prove this by contradiction. Suppose $G$ is in one of these
cases. By replacing $M$ by a two-fold cover in case~\ref{O} (the
resulting $M$ also satisfies the hypothesis as can be seen by pulling
back contact structures on $M$), we can assume that we are in one of
the cases~\ref{T}, \ref{I} and \ref{T2}. As $M$ has a universally
tight positive contact $\xi$, as before after passing to a cover $N$
(with $\abs{\pi_1(N)}\geq 4$) the pullback of $\xi$ has trivial Euler class
and hence a framing  associated to it. This pulls back to a
framing isotopic to left-invariant framing on $S^3$ as $\xi$ lifts to
the standard contact structure on $S^3$. Thus, $N$ has a framing that
pulls back to the left-invariant framing of $S^3$

But, as $-M$ also has a universally tight positive contact structure
$\xi'$, we can apply the same argument to $-M$. Thus, the pullback of
$\xi$ to $-N$ has associated with it a framing, and this pulls back to
a framing isotopic to the left-invariant framing. By the above
proposition, it follows that $N$ has a framing that pulls back to the
right-invariant framing of $S^3$. This contradicts
Corollary~\ref{T:lftrt}
\end{proof}

\end{document}